# Correction to a lemma in
## *Embedded contact homology and Seiberg-Witten Floer homology* IV

Clifford Henry Taubes[†]
Department of Mathematics
Harvard University

chtaubes@math.harvard.edu

ABSTRACT: This note corrects an erroneous statement in Lemma 3.8 of the author's paper *Embedded contact homology and Seiberg-Witten Floer homology IV* which was published in Volume 14 of Geometry and Topology in 2009.

[†] Supported in part by a grant from the National Science Foundation

# 1. Introduction

The purpose of this note is to correct a mistatement in a lemma from [1] that enters the proof of the equivalence, for a 3-manifold with a contact structure, between the Seiberg-Witten Floer cohomology and Michael Hutching's embedded contact homology. This is Lemma 3.8 in [1]. A precise statement of the correction is given below in Sections 1b. What follows directly in Section 1a sets background and notation.

## a) Background and notation

Let M denote a compact, 3-dimensional manifold with a given contact 1-form. This 1-form is denoted by $a$ in what follows. To say that $a$ is a contact form is to say that $a \wedge da \neq 0$ at all points of M. The convention is to orient M using the 3-form $a \wedge da$. A Riemannian metric on M is chosen so that $a$ has norm 1 and so that its Hodge star is $\frac{1}{2} da$.

Fix a $\mathrm{Spin}^{\mathbb{C}}$ structure for M which is a lift of its oriented orthonormal frame bundle to a principal U(2) bundle (which is $\mathrm{Spin}^{\mathbb{C}}(3)$). Let $\mathbb{S} \to M$ denote the associated $\mathbb{C}^2$ bundle to this principle U(2) bundle. This is a module for Clifford multiplication which is a homomorphism to be denoted by cl from T*M to vector space of anti-Hermitian endomorphisms of $\mathbb{S}$. The convention in this paper (and in [1]) has cl obeying the rule

$$\mathrm{cl}(u)\,\mathrm{cl}(v) = -\langle u, v \rangle - \mathrm{cl}(*(u \wedge v)).$$

(1.1)

(Note the minus sign in the right most term.)

The covariant derivative acting on sections of $\mathbb{S}$ can be defined from the Levi-Civita connection on M's orthonormal frame bundle and a given Hermitian connection on the complex line bundle $\det(\mathbb{S})$ (which is $\wedge^2 \mathbb{S}$, which is associated to the principal U(2) bundle via the determinant homomorphism from U(2) to U(1)). Supposing that $\mathbb{A}$ is a connection on $\det(\mathbb{S})$, then the corresponding covariant derivative, $\nabla^{\mathbb{A}}$, maps sections of $\mathbb{S}$ to sections of $\mathbb{S} \otimes T^*M$. Composing this with the Clifford multiplication map gives a Dirac operator, $D^{\mathbb{A}}$, which is a first order, differential operator mapping sections of $\mathbb{S}$ to sections of $\mathbb{S}$.

A pair $(\mathbb{A}, \Psi)$ of Hermitian connection on $\det(\mathbb{S})$ and section of $\mathbb{S}$ is said to obey the Seiberg-Witten equations on M when

$$B_{\mathbb{A}} = \Psi^{\dagger}\tau\Psi \quad \textit{and} \quad D^{\mathbb{A}}\Psi = 0$$

(1.2)

where $B_{\mathbb{A}}$ denotes the metric Hodge dual of the curvature 2-form of $\mathbb{A}$ (the former is an $i\mathbb{R}$ valued 1-form because the latter is an $i\mathbb{R}$ valued 2-form) and where $\Psi^{\dagger}\tau\Psi$ is shorthand for the $i\mathbb{R}$ valued 1 form that is defined by the following rule: Its inner product with any given 1-form $u$ is equal to the inner product of $\Psi$ with $cl(u)\Psi$.

A variant of the Seiberg-Witten equations is defined with the choice of a coclosed, $i\mathbb{R}$ valued 1-form on M. (A 1-form $\varpi$ is coclosed when $d*\varpi = 0$.) Supposing that $\varpi$ is such a form, the variant has $(\mathbb{A},\psi)$ obeying

$$B_{\mathbb{A}} = \Psi^{\dagger}\tau\Psi + \varpi \quad and \quad D^{\mathbb{A}}\Psi = 0$$

(1.3)

Of interest in what follows are the versions of (1.3) for a 1-parameter family of choices for $\varpi$ that is defined as follows: The family is parametrized by $[1,\infty)$; and any given $r \in [1,\infty)$ member of the family has $\varpi = 2(-ira + *d\mu)$ with $\mu$ denoting some $r$-independent, $i\mathbb{R}$ valued 1-form and $*$ denoting the Hodge dual from the metric on M.

To continue introducing the background notions, note that multiplication by the 1-form $a$ defines the splitting of $\mathbb{S}$ as an orthogonal, direct sum of line bundles:

$$\mathbb{S} = E \oplus EK^{-1}.$$

(1.2)

Here, the E summand is the $+i$ eigenbundle for $cl(a)$; and $EK^{-1}$ is the $-i$ eigenbundle. (The complex line bundle $K^{-1}$ is isomorphic as an oriented 2-dimension real vector bundle to the kernel of $a$ in TM.) A section of $\mathbb{S}$ will be written at times as a pair $(\alpha, \beta)$ with $\alpha$ being the part in the E summand of (1.2) and with $\beta$ being the part in the $EK^{-1}$ summand.

The canonical $\text{Spin}^{\mathbb{C}}$ structure on M, as defined by $a$, is the $\text{Spin}^{\mathbb{C}}$ structure that gives (1.2) with E having zero first Chern class. This is to say that E is isomorphic to the product complex line bundle (which is denoted by $\underline{\mathbb{C}}$ in what follows). The corresponding version of $\mathbb{S}$ is denoted by $\mathbb{S}_I$. Let $1_{\mathbb{C}}$ denote a chosen, unit length section of $\underline{\mathbb{C}}$. This defines a section of $\mathbb{S}_I$ (to be denoted by $\Psi_I$) that is written using the $E = \underline{\mathbb{C}}$ version of (1.2) as $(1_{\mathbb{C}}, 0)$. As explained in [2], there is a unique connection on $det(\mathbb{S}_I) = K^{-1}$ with the property that $\Psi_I$ obeys the corresponding Dirac equation. This connection is denoted by $\mathbb{A}_K$. Now, supposing that $\mathbb{S}$ comes from *some* $\text{Spin}^{\mathbb{C}}$ structure, then any connection on $det(\mathbb{S})$ can be written unambiguously as $\mathbb{A}_K + 2A$ where A is a connection on the corresponding version of the bundle E. The corresponding version of the Dirac operator $D^{\mathbb{A}}$ is denoted in what follows by $D_A$.

It proves useful to write the Seiberg-Witten equations in (1.3) with $\varpi$ given by the 1-form $2(-ira + *d\mu)$ as equations for the pair $(A,\psi)$ with $\psi = (2r)^{1/2}\Psi$:

$$B_A = r((\psi^\dagger \tau \psi) - i a) + *d\mu - \tfrac{1}{2} B_{\mathbb{A}_K} \quad \text{and} \quad D_A \psi = 0 \;,$$

(1.3)

This is the form of the Seiberg-Witten equations on M that is considered in [1].

The particular lemma in [1] (which is Lemma 3.8) concerns not so much solutions to (1.3) (although it applies to them), but solutions to the analog of (1.3) for an $\mathbb{R}$-parametrized family $\{(A,\psi)|_s: s \in \mathbb{R}\}$ of connection on E and section of $\mathbb{S}$. The analog asks that this 1-parameter family of connection on E and section of $\mathbb{S}$ obey

$$\partial_t A + *B_A = r((\psi^\dagger \tau \psi) - ia) + *d\mu - \tfrac{1}{2} B_{\mathbb{A}_K} \quad \text{and} \quad \partial_t \psi + D_A \psi = 0$$

(1.4)

at each $s \in \mathbb{R}$ and have $s \to \infty$ and $s \to -\infty$ limits that obey (1.3). An *instanton* solution to (1.4) is just such a smoothly parametrized family of pairs, $\{(A, \psi)|_s: s \in \mathbb{R} \in \mathbb{R}\}$ of connection on E and section of $\mathbb{S}$ that obeys (1.4) with $s \to \infty$ and $s \to -\infty$ limits that obey (1.3). The notation in what follows uses $\mathfrak{d}$ to denote any given instanton solution to (1.4).

**b) The correction to Lemma 3.8 in [1]**

There are two real numbers that are associated to any given instanton solution to (1.4); if $\mathfrak{d}$ is the instanton, they are denoted by $A_\mathfrak{d}$ and $F_\mathfrak{d}$. Both are defined in Section 3.1 of [1]. Their precise definitions are not of concern in what follows except to the extent that they are referred to in Lemma 3.8 in [1], and they are referred to there only because Lemmas 3.2-3.7 in [1] can be invoked when certain bounds on $A_\mathfrak{d}$ or $F_\mathfrak{d}$ are assumed. Lemma 3.8 in [1] also refers to a function on $\mathbb{R}$ that is defined from an instanton solution to (1.4). The function is denoted by $\underline{M}$ and it is defined as follows: Let $\mathfrak{d}$ denote the instanton in question. Write the corresponding section $\psi$ using (1.2) as a pair $(\alpha, \beta)$ with $\alpha$ (as always) being the part in the E summand. For $s \in \mathbb{R}$,

$$\underline{M}(s) = r \int_{[s-1,s-1] \times M} (1 - |\alpha|^2) \;.$$

(1.5)

By way of an explanation: Lemma 3.1 in [1] says in part that the function $|\alpha|^2$ obeys $|\alpha|^2 \leq 1 + \mathcal{O}(\tfrac{1}{r})$. As a consequence, $\underline{M}$ has an a priori $r$-independent lower bound. Therefore, any upper bound (less than $\mathcal{O}(r)$) for $\underline{M}(s)$ constrains the size of the subset in $[s-1, s+1] \times M$ where $|\alpha|$ is significantly less than 1.

By way of one last bit of notation: Supposing that A is now a Hermitian connection over $\mathbb{R} \times M$ on the bundle E (the bundle E is pulled back by the projection map to M), and $\psi = (\alpha, \beta)$ is a section over $\mathbb{R} \times M$ of the bundle $\mathbb{S}$ (also pulled back by the projection map), let $\nabla_A \alpha$ and $\nabla_A \beta$ now denote the respective covariant derivatives of $\alpha$ and $\beta$ on $\mathbb{R} \times M$ as defined by the connection A on E in the case of $\alpha$, and by the connection A and the pull-back of $A_K$ in the case of $\beta$. To be sure: These covariant derivatives involve derivatives along the $\mathbb{R}$ factor of $\mathbb{R} \times M$ as well as the M factor. Note in particular that the 1-parameter family of connections from an instanton solution to (1.4) canonically defines a connection on E over $\mathbb{R} \times M$; and the 1-parameter family of sections of $\mathbb{S}$ from an instanton defines a section of $\mathbb{S}$ over $\mathbb{R} \times M$.

What follows directly is the *corrected* statement of Lemma 3.8 in [1]. The misstated part of the lemma in [1] is given subsequently in (1.6).

**Lemma 3.8** (CORRECTED): *Given $\mathcal{K} \geq 1$, there exists $\kappa \geq 1$ with the following significance: Suppose that $r \geq \kappa$, and that $\mathfrak{d} = (A, \psi = (\alpha, \beta))$ is an instanton solution to (1.4) with either $A_\mathfrak{d} < r^2$ or $F_\mathfrak{d} \geq -r^2$. Fix a point $s_0 \in \mathbb{R}$ and $R \geq 2$; and suppose that $\sup_{s \in [s_0 - R - 3, s_0 + R + 3]} \underline{M}(s) \leq \mathcal{K}$. Let $X_* \subset \mathbb{R} \times M$ denote the subset of points where $1 - |\alpha| \geq \kappa^{-1}$. The bounds stated below hold on the domain $[s_0 - R, s_0 + R] \times M$.*

- $|\nabla_A \alpha|^2 + r|\nabla_A \beta|^2 \leq \kappa r(1 - |\alpha|^2) + \kappa^2$
- $(1 - |\alpha|^2) \leq \kappa(\frac{1}{r} + e^{-\sqrt{r}\, \text{dist}(\cdot, X_*)/\kappa})$
- $|\nabla_A \alpha|^2 + r|\nabla_A \beta|^2 \leq \kappa(\frac{1}{r} + r\, e^{-\sqrt{r}\, \text{dist}(\cdot, X_*)/\kappa})$
- $|\beta|^2 \leq \kappa \frac{1}{r}(\frac{1}{r} + e^{-\sqrt{r}\, \text{dist}(\cdot, X_*)/\kappa})$

The version of Lemma 3.8 in [1] has the second bullet stating that

$$r(1 - |\alpha|^2) \leq \kappa(\frac{1}{r} + r\, e^{-\sqrt{r}\, \text{dist}(\cdot, X_*)/\kappa}) \,.$$

(1.6)

Note the extra factor of $r$ on the left. This is the only misstatement in the lemma (to the author's knowledge.) The corrected version of the lemma is proved in the upcoming Section 2 of this note. With regards to (1.6), one other remark is in order: The upper bound in (1.6) does infact hold when $\mathfrak{d} = (A, \psi)$ obeys the modified version of (1.4) with both $\mu$ and $B_{\mathbb{A}_K}$ absent from the left most equation.

The original version of Lemma 3.8 is used in other sections of [1]. Additional remarks are needed with regards to these applications when using the corrected version of Lemma 3.8. There remarks are given momentarily as they require some additional stage setting. To this end, let $\wp$ denote a smooth, non-negative and non-decreasing function on $[0, \infty)$ which obeys $\wp(x) = x$ for $x$ near 0 and $\wp(1) = 1$. Supposing that A is a connection on E's pull-back over $\mathbb{R} \times M$ and $\alpha$ is a section over $\mathbb{R} \times M$ of this pull-back, use $\wp$ to define a new connection (denoted by $\hat{A}$):

$$\hat{A} = A - \tfrac{1}{2}\wp(|\alpha|^2)|\alpha|^{-2}(\bar{\alpha}\nabla_A\alpha - \alpha\nabla_A\bar{\alpha}),$$

(1.7)

with the covariant derivative $\nabla_A$ here and subsequently denoting the covariant derivative on $\mathbb{R} \times M$ as defined by the connection A. (To be sure: The covariant derivative $\nabla_A$ differentiates not just along the tangents to the $s \in \mathbb{R}$ slices, but along the $\mathbb{R}$ direction.) The curvature 2-form of $\hat{A}$ is given by the formula

$$F_{\hat{A}} = (1 - \wp)F_A - \wp'\nabla_A\bar{\alpha} \wedge \nabla_A\alpha.$$

(1.8)

The connection $\hat{A}$ is introduced in Equation (3.33) of [1]; and equation (3.35) in [1] makes the following assertion:

*If the assumptions in Lemma 3.8 hold, then* $|F_{\hat{A}}| \leq c_0(\tfrac{1}{r} + re^{-\sqrt{r}\,\mathrm{dist}(\cdot,X_*)/\kappa})$
*at all points with* $s \in [s_0 - R, s_0 + R]$.

(1.9)

Here and in [1], $c_0$ denotes a number greater than 1 that is independent of $r$ and the instanton $(A, \psi)$. (The verbatim statement in [1] says that (1.9) holds if the assumptions of Lemma 3.<u>6</u> hold. This is a typo. It should have been Lemma 3.8.)

The subsequent appeals to the original version of Lemma 3.8 use either the bound in (1.9) or they follow using the bounds given above in the corrected version of Lemma 3.8. Meanwhile, (1.9) follows using the second and third bullets in the corrected version of Lemma 3.8 above with the assertion made by Lemma 3.9 in [1]. (The proofs of Lemmas 3.9 and 3.10 in [1] makes no appeal to the erroneous bound in (1.6).)

## 2. Proof of the corrected version of Lemma 3.8

The proof is given in the three subsections that follow. (A version of the proof of the correction is also given in the Appendix to [3].) The notation in what follows has $c_0$ denoting a number that is greater than 2 and independent of the particular instanton $\mathfrak{d}$ and independent of $r$. The particular value of $c_0$ can be assumed to increase between successive appearances.

### a) Initial remarks

The first bullet in the top set of four bullets is proved using the argument for the proof of the corresponding first bullet in Lemma 3.8 of [1]. (That argument in [1] is fine.) The second bullet of Lemma 3.8 in [1] makes the assertion that

$$r(1 - |\alpha|^2) + |\nabla_A \alpha|^2 + r|\nabla_A \beta|^2 \leq \kappa(\tfrac{1}{r} + r e^{-\sqrt{r}\, \text{dist}(\cdot, X_*)/\kappa})$$

(2.1)

which is the assertion of the third bullet of the corrected Lemma 3.8 given above plus the erroneous assertion in (1.6). The argument in [1] that claims to prove (2.1) does in fact prove the weaker assertion:

$$r(1 - |\alpha|^2) + |\nabla_A \alpha|^2 + r|\nabla_A \beta|^2 \leq \kappa(1 + r e^{-\sqrt{r}\, \text{dist}(\cdot, X_*)/\kappa}).$$

(2.2)

This implies what is asserted by the second bullet of the corrected Lemma 3.8. But, it does not imply by itself what is asserted by the third bullet.

**b) The bound for $|\nabla_A \alpha|$**

The argument for the bound

$$|\nabla_A \alpha|^2 \leq \kappa(\tfrac{1}{r} + r e^{-\sqrt{r}\, \text{dist}(\cdot, X_*)/\kappa})$$

(2.3)

is given directly in seven steps.

Step 1: Let $\kappa_*$ denote the version of the number $\kappa$ that appears in the first and second bullets of corrected Lemma 3.8. The bound that is asserted by (2.3) holds (with $\kappa = c_0 \kappa_*$) where the distance to $X_*$ is less than $100\kappa_* \sqrt{r}$ because of the first two bullets in the lemma.

Step 2: According to the second bullet of the corrected Lemma 3.8, the norm of $\alpha$ obeys $|\alpha| \geq \tfrac{99}{100}$ where the distance to $X_*$ is greater than $c_0 \kappa_* \sqrt{r}$ when $r > c_0$ (which will henceforth be assumed). Meanwhile, the bundle E where $|\alpha| > \tfrac{99}{100}$ is isomorphic to the product $\mathbb{C}$ bundle by an isomorphism that takes $\alpha$ to a real number that can be written as $1 - z$ with $z < \tfrac{1}{100}$. By virtue of the second bullet of the corrected Lemma 3.8 and Lemma 3.1 in [1], this function $z$ obeys

$$|z| \leq c_0(\tfrac{1}{r} + e^{-\sqrt{r}\, \text{dist}(\cdot, X_*)/c_0}).$$

(2.4)

(Lemma 3.1 in [1] asserts in part the bound $-z \leq c_0 \tfrac{1}{r}$.)

The isomorphism from E to the product $\mathbb{C}$ bundle that depicts $\alpha$ as $(1 - z)$ identifies $\beta$ with a section of $K^{-1}$ to be denoted by $\beta_\Diamond$. This isomorphism also identifies A with a connection that is written as $\theta_0 + \hat{a}_A$ with $\theta_0$ denoting the product connection on the product $\mathbb{C}$ bundle and with $\hat{a}_A$ denoting an $i\mathbb{R}$-valued 1-form.

This same isomorphism identifies $\nabla_A\alpha$ with a $\mathbb{C}$-valued 1-form that can be written as $-dz + \hat{a}_A(1-z)$. Thus, bounds for $|dz|$ and $|\hat{a}_A|$ lead to bounds for $|\nabla_A\alpha|$ (and vice-versa because $z$ is real and $\hat{a}_A$ is $i\mathbb{R}$ valued). In particular, a $c_0 \frac{1}{\sqrt{r}}(1 + r e^{-\sqrt{r}\, \text{dist}(\cdot, X_*)/c_0})$ bound for the norms of $dz$ and $\hat{a}_A$ implies what is asserted by (2.3).

Step 3: The Riemannian metric on $\mathbb{R} \times M$ with the self-dual 2-form $ds \wedge a + *a$ defines an almost complex structure on $T(\mathbb{R} \times M)$ which is $\mathbb{R}$-invariant and maps $\frac{\partial}{\partial s}$ to the Reeb vector field. (The $\mathbb{R}$ action is translation along the $\mathbb{R}$ factor.) With this almost complex structure understood, introduce by way of notation $\bar{\partial}_A\alpha$ to denote $T^{0,1}$ part of $\nabla_A\alpha$.

The right most equation in (1.4) when written in terms of $\alpha$ and $\beta$ identifies $\bar{\partial}_A\alpha$ with $x(\nabla_A\beta)$ where $x$ is a homormorphism with norm bounded by $c_0$. Using the isomorphism of the preceding step, this identification takes the form

$$\bar{\partial} z + \hat{a}_A^{0,1}(1-z) = x(\nabla_A\beta) .$$

(2.5)

with $\bar{\partial} z$ denoting the $T^{0,1}$ part of $dz$ and with $\hat{a}_A^{0,1}$ denoting the analogous part of $\hat{a}_A$. Because $\hat{a}_A$ is $i\mathbb{R}$ valued, the norm of $\hat{a}_A^{0,1}$ bounds the norm of $\hat{a}_A$. Therefore, the identity in (2.5) and the inequality in (2.2) lead to the bound

$$|\hat{a}_A| \leq c_0 |\bar{\partial} z| + c_0 \frac{1}{\sqrt{r}}(1 + \sqrt{r}\, e^{-\sqrt{r}\, \text{dist}(\cdot, X_*)/c_0}) .$$

(2.6)

This implies that a bound on $|dz|$ by $c_0 \frac{1}{\sqrt{r}}(1 + r e^{-\sqrt{r}\, \text{dist}(\cdot, X_*)/c_0})$ leads to the same bound on $|\hat{a}_A|$ (perhaps with a large $c_0$). Steps 4-7 derive this desired bound on the norm of $dz$.

Step 4: The right most equation in (1.4) says that $(\frac{\partial}{\partial s} + D_A)\psi = 0$, and thus that

$$(-\tfrac{\partial}{\partial s} + D_A)(\tfrac{\partial}{\partial s} + D_A)\psi = 0 .$$

(2.7)

This last identity, when written in terms of $\alpha$ and $\beta$ and projected to the respective $E$ and $EK^{-1}$ summands of the spinor bundle has the form

- $\nabla_A^\dagger \nabla_A \alpha + r(|\alpha|^2 - 1 + |\beta|^2)\alpha + c_0\alpha + c_1\nabla_A\beta + c_2\beta = 0$,
- $\nabla_A^\dagger \nabla_A \beta + r(|\alpha|^2 + 1 + |\beta|^2)\beta + c_3\nabla_A\beta + c_4\beta + c_5\nabla_A\alpha + c_6\alpha = 0$,

(2.8)

where $\{c_k\}_{k=0,\ldots,6}$ are endomorphism that are independent of $(A,\psi)$ and $r$ and have norms bounded by $c_0$. (By way of a parenthetical remark, the endomorphisms $c_0$ and $c_6$ have linear dependence on $*d\mu$ and $B_{A_K}$.)

Because of the top bullet of (2.8), the function $z$ obeys an equation that can be written schematically as

$$d^\dagger dz + 2rz = -|\hat{a}_A|^2(1-z) + r((1-|\alpha|^2 + |\beta|^2 + z)z + x_0(1-z) + x_1(\nabla_A\beta) + x_2\beta \tag{2.9}$$

with $x_0$, $x_1$, and $x_2$ each denoting a homomorphism whose norm is bounded by $c_0$. (Equation (2.9) comes from taking the real part of (2.8) after E is identified with the product bundle in the manner described previously.)

Step 5: With (2.9) in hand, fix a point $p \in \mathbb{R} \times M$ with distance at least $c_0 \kappa \sqrt{r}$ from the complement of $X_*$. Let $\rho$ denote this distance to the complement of $X_*$, and let B denote the ball of radius equal to the minimum of $\frac{1}{2}\rho$ and $c_0^{-1}$ centered at the point p.

Define the bump function $\chi_p$ by the rule $\chi_p(\cdot) = \chi(4\rho^{-1}\mathrm{dist}(\cdot,p))$. This function is equal to 1 where the distance to p is less than $\frac{1}{16}\rho$ and it is equal to zero where the distance to p is greater than $\frac{1}{4}\rho$. Let $z_p = \chi_p z$. This function has compact support in B; and, by virtue of (2.9), it obeys an equation that has the schematic form

$$d^\dagger dz_p + 2rz_p = -2\langle d\chi_p, dz\rangle + d^\dagger d\chi_p z + \chi_p h \tag{2.10}$$

where $\langle\,,\rangle$ denotes the Riemannian inner product on cotangent vectors; and where $h$ is shorthand for what appears on the right hand side of (2.9).

Step 6: Let q denote for the moment a point in B with distance less than $\frac{1}{32}\rho$ from p. Introduce now $G_q$ to denote the Dirichelet Green's function for the operator $d^\dagger d + \frac{1}{2}r$ on B with pole at the point q. This Green's function is zero on the boundary of B, positive inside B and smooth except at p. Moreover, it obeys:

- $G_q \leq c_0 \frac{1}{\mathrm{dist}(\cdot,q)^2} e^{-\sqrt{r}\,\mathrm{dist}(\cdot,q)/c_0}$
- $|dG_q| + |\nabla_q G_q| \leq c_0 \frac{1}{\mathrm{dist}(\cdot,q)^3} e^{-\sqrt{r}\,\mathrm{dist}(\cdot,q)/c_0}$ .
- $|(\nabla_q(dG_q)| \leq c_0 \frac{1}{\mathrm{dist}(\cdot,q)^4} e^{-\sqrt{r}\,\mathrm{dist}(\cdot,q)/c_0}$ .

(2.11)

Here, $\nabla_q$ denotes differentiation with respect to the pole position, q. Meanwhile, dG denotes differentiation of G with respect to the argument of the function $G_q(\cdot)$. The last two bounds follow using the bound $xe^{-x} \leq c_0 e^{-x/2}$ (which holds for any $x \geq 0$) to replace factors of $\sqrt{r}$ with inverse powers of $dist(\cdot, q)$.

Step 7: The Green's function $G_q$ is used with (2.10) to write $z_p$ at any point q in the radius $\frac{1}{256}\rho$ ball centered at p as

$$z_p|_q = \int_B G_q(-2\langle dz, d\chi_p \rangle + d^\dagger d\chi_p z + \chi_p h)$$

(2.12)

Keeping in mind that $d\chi_p$ is non-zero only where the distance to q is at least $\frac{1}{16}\rho$, an integration by parts can be used to remove the derivative from $z$ that appears in the integral on the left hand side of (2.12) with the result being:

$$z_p|_q = \int_B (2\langle dG_q, d\chi_p \rangle - G_q d^\dagger d\chi_p)z + G_q \chi_p h).$$

(2.13)

The derivative of this identity (with respect to q) gives an identity for $dz|_q$ that has the same form but with $\nabla_q G_q$ replacing $G_q$. As explained in the next two paragraphs, the desired bound

$$|dz|_p| \leq c_0 \frac{1}{\sqrt{r}}(1 + re^{-\sqrt{r}\,dist(\cdot, X_*)/c_0})$$

(2.14)

then follows from the $q = p$ version of the latter identity using (2.11) with the inequalities in (2.4) and (2.2) and the bound $|\beta| \leq c_0 \frac{1}{\sqrt{r}}$ from Lemma 3.1 in [1].

To see in more detail how the bound in (2.14) comes about, first consider bounding the norm of the contribution to $dz|_p$ from the evaluation at $q = p$ of the term

$$\int_B (2\langle \nabla_q dG_q, d\chi_p \rangle - \nabla_q G_q d^\dagger d\chi_p)z.$$

(2.15)

By virtue of (2.11), both $|\langle \nabla_q dG_q, d\chi_p \rangle|$ and $|\nabla_q G_q d^\dagger d\chi_p|$ at $q = p$ are bounded by

$$c_0 \frac{1}{\rho^5} e^{-\sqrt{r}\,\rho/c_0}.$$

(2.16)

Meanwhile, $|z| \leq c_0(\frac{1}{r} + e^{-\sqrt{r}\,\rho/c_0})$ by virtue of (2.4). Therefore, since the volume of B is at most $c_0\rho^4$, the norm of (2.15) is at most $c_0 \frac{1}{r\rho} e^{-\sqrt{r}\,\rho/c_0}$ which is consistent with (2.17) since $\rho \geq \frac{1}{\sqrt{r}}$.

Consider now the contribution to $|dz_p|$ from the integral of $\nabla_q G_q \chi_\rho h$. The inequalities in (2.2) and (2.4) lead to the bound $|h| \leq c_0(1 + r\, e^{-\sqrt{r}\, \rho/c_0})$ on the support of $\chi_\rho$. It follows from this bound that the norm of the integral over B of $\nabla_q G_q \chi_\rho h$ is at most $c_0(1 + r\, e^{-\sqrt{r}\, \rho/c_0})$ times the integral of $|\nabla_q G_q|$ over the support of $\chi_\rho$. And, the latter integral is at most $c_0 \frac{1}{\sqrt{r}}$ by virtue of the second bullet in (2.11).

### c) Bounds for $|\beta|$ and $|\nabla_A \beta|$

This subsection derives the bounds for $|\beta|$ and for $|\nabla_A \beta|$ that are asserted by the third and fourth bullets of the corrected version of Lemma 3.8. To this end: The $|\beta|^2$ bound follows directly from the bound in (2.2) for $(1 - |\alpha|^2)$ and what is said by the second bullet of Lemma 3.1 in [1]. (This bullet says that $|\beta|^2 \leq c_0 \frac{1}{r}(1 - |\alpha|^2) + c_0 \frac{1}{r^2}$.)

The proof of the bound for $|\nabla_A \beta|^2$ uses the equation in the second bullet of (2.8) with a Green's function that obeys the bounds in (2.11). To explain in detail, note first that $|\nabla_A \beta|$ is the same as $|\nabla_{\theta_0 + \hat{a}_A} \beta_\diamond|$ at points where the distance to $X_*$ is greater than $c_0 \frac{1}{\sqrt{r}}$. Let p denote such a point and let B denote the ball centered at p from Step 5 in the previous subsection. If the radius of B is less than $c_0^{-1}$ (which can be assumed), then the bundle $K^{-1}$ on B is isomorphic to the product $\mathbb{C}$ bundle. Moreover, one can and should fix an isomorphism that identifies the canonical connection on $K^{-1}$ with a connection on the product bundle that can be written as $\theta_0 + \Gamma$ with $|\Gamma|$ and $|\nabla \Gamma|$ both less than $c_0$. Having done this, then writing the second bullet of (2.8) using this isomorphism gives and equation for $\beta_\diamond$ that has the schematic form

$$\nabla_{\theta_0}^\dagger \nabla_{\theta_0} \beta_\diamond + r\, \beta_\diamond = c_6 + \mathcal{R}$$

(2.17)

where $c_6$ is the same as its namesake in (2.8) after accounting for the isomorphism between $K^{-1}$ and the product bundle; and where $\mathcal{R}$ obeys the pointwise bound

$$|\mathcal{R}| \leq c_0 \frac{1}{\sqrt{r}} (1 + r\, e^{-\sqrt{r}\, \text{dist}(\cdot,\, X_*)/c_0}).$$

(2.18)

With regards to $c_6$: It appears when the term $c_6 \alpha$ in the second bullet of (2.18) is written as $c_6(1 - z)$. The $-c_6 z$ part of the latter is encorporated into the $\mathcal{R}$ term in (2.17). With regards to the norm of $\mathcal{R}$: The only contribution to $\mathcal{R}$ that does not directly obey the bound in (2.18) by virtue of what has been proved up to this point is a term $-(d^\dagger \hat{a}_A) \beta_\diamond$ that comes from writing out $(\nabla_{\theta_0 + \hat{a}_A})^\dagger \nabla_{\theta_0 + \hat{a}_A} \beta_\diamond$. To see about the size of $|(d^\dagger \hat{a}_A) \beta_\diamond|$ note the following: Whereas the real part of the top bullet in (2.8) (after E is identified with the product $\mathbb{C}$ bundle as described previously) gives (2.9), the corresponding imaginary part of (2.9) asserts that

$$d^\dagger \hat{a}_A (1-z) - 2\langle \hat{a}_A, dz \rangle + \mathrm{im}(c_1 \nabla_{\theta_0 + \hat{a}_A} \beta_\diamond + c_2 \beta) = 0 .$$

(2.19)

This identity with (2.2) and with the previously derived bounds on $|\hat{a}_A|$ and $|dz|$ leads to a bound for $|(d^\dagger \hat{a}_A)\beta_\diamond|$ that is consistent with the asserted bound for $|\mathcal{R}|$.

With (2.17) in hand, then arguments much like those used to bound $dz$ in the last steps of the previous subsection (using a Green's function obeying (2.11) and the function $\chi_\rho$), can be used with only one extra remark to bound the norm of $\nabla_{\theta_0} \beta_\diamond$:

$$|\nabla_{\theta_0} \beta_\diamond| \leq c_0 \tfrac{1}{r} (1 + r e^{-\sqrt{r}\, \mathrm{dist}(\cdot, X_*)/c_0}).$$

(2.20)

The extra remark is this: The solution to the equation $\nabla_{\theta_0}^\dagger \nabla_{\theta_0} w + \tfrac{1}{2} rw = c_6$ that vanishes on $\partial B$ is such that both its norm *and* the norms of its first derivatives are pointwise bounded at p by $c_0 \tfrac{1}{r}$.

With (2.20) in hand, then the inequality $|\nabla_A \beta| \leq |\nabla_{\theta_0} \beta_\diamond| + |\hat{a}_A||\beta_\diamond|$ with the bounds proved previously for $|\hat{a}_A|$ and $|\beta|$ lead directly to the desired bound for the norm of $\nabla_A \beta$.